\documentclass{amsart}
\usepackage{amssymb}
\usepackage{amsfonts}
\usepackage{latexsym}
\newcommand{\g}{{\mathfrak g}}
\newcommand{\h}{{\mathfrak h}}
\newcommand{\n}{{\mathfrak n}}

\newtheorem{theorem}{Theorem}[section]
\newtheorem{lemma}[theorem]{Lemma}

\theoremstyle{definition}

\newtheorem{example}[theorem]{Example}

\newtheorem{conjecture}[theorem]{Conjecture}

\begin{document}
\author{Pavel Etingof, Massachusetts Institute of Technology,
  USA}
\title{On the dynamical Yang-Baxter equation}
\maketitle
\centerline{\bf To Mira}
\vskip .05in
\centerline{\bf 2000 Mathematics Subject Classification: 17B37,20G42,81R50}

\section{Introduction}

This talk is inspired by two previous ICM talks, \cite{Dr} and
\cite{Fe}. 
Namely, one of the main ideas of
\cite{Dr} is that the quantum Yang-Baxter equation
(QYBE), which is an important
equation arising in quantum field theory and statistical
mechanics, is best understood within the framework of Hopf
algebras, or quantum groups. On the other hand, in \cite{Fe},
it is explained that another important equation of mathematical 
physics, the star-triangle relation, may (and should) be viewed
as a generalization of QYBE, in which
solutions depend on additional ``dynamical'' parameters. It is 
also explained there that to a solution of the quantum 
dynamical Yang-Baxter equation one may associate a kind of
quantum group. These ideas gave rise to a vibrant new 
branch of ``quantum algebra'', which may be called 
the theory of dynamical quantum groups. My goal in this talk 
is to give a bird's eye review of some aspects of this theory and its
applications.   

{\bf Remark.} Because of size restrictions, the list of references 
below is not complete. The reader is referred to original papers 
and review \cite{ES2} for more references. 

{\bf Acknowledgements.} I am indebted to G. Felder 
for creating this subject, and deeply grateful to I. Frenkel, A. Kirillov Jr., 
A. Moura, O. Schiffmann, and A. Varchenko for 
collaboration that led to this work. 
This work was partially supported by the NSF grant DMS-9988796, 
and was done in part for the Clay Mathematics Institute. 
I am grateful to the Max Planck Institute for hospitality. 

\section{The quantum dynamical Yang-Baxter equation}

\subsection{The equation}
The quantum dynamical Yang-Baxter equation (QDYBE), proposed by Gervais, 
Neveu, and Felder, is an equation 
with respect to a (meromorphic) function $R:{\mathfrak h^*}\to {\rm
  End}_{\mathfrak h}(V\otimes V)$, where ${\mathfrak h}$ is a
commutative finite dimensional Lie algebra over $\mathbb C$, and $V$ is a 
semisimple ${\mathfrak h}$-module. It reads
$$
R^{12}(\lambda-h^3)R^{13}(\lambda)R^{23}(\lambda-h^1)=
R^{23}(\lambda)R^{13}(\lambda-h^2)R^{12}(\lambda)
$$
on $V\otimes V\otimes V$. Here $h^i$ is the {\it dynamical notation}, 
to be extensively used below: for instance, $R^{12}(\lambda-h^3)$ is defined
by the formula $R^{12}(\lambda-h^3)(v_1\otimes v_2\otimes v_3):=
\bigl(R^{12}(\lambda-\mu)(v_1\otimes v_2)\bigr)\otimes v_3$
if $v_3$ is of weight $\mu$ under $\mathfrak h$. 

It is also useful to consider 
QDYBE with spectral parameter. 
In this case the unknown function $R$ depends meromorphically 
on an additional complex variable $u$. Let $u_{ij}=u_i-u_j$. Then 
the equation reads
$$
R^{12}(u_{12},\lambda-h^3)R^{13}(u_{13},\lambda)R^{23}(u_{23},\lambda-h^1)=
R^{23}(u_{23},\lambda)R^{13}(u_{13},\lambda-h^2)R^{12}(u_{12},\lambda).
$$

Sometimes it is necessary to consider QDYBE with step $\gamma\in \mathbb C^*$, 
which differs from the usual QDYBE by the replacement 
$h^i\to \gamma h^i$. Clearly, $R(\lambda)$ satisfies QDYBE 
iff $R(\lambda/\gamma)$ satisfies QDYBE with step $\gamma$. 

Invertible solutions of QDYBE are called quantum dynamical 
R-matrices (with or without spectral parameter, and with step 
$\gamma$ if needed). If ${\mathfrak h}=0$, 
QDYBE turns into the usual quantum Yang-Baxter equation
$R^{12}R^{13}R^{23}=R^{23}R^{13}R^{12}$. 

\subsection{Examples of solutions of QDYBE}
Let
$V$ be the vector representation of $sl(n)$, and $\h$ 
the Lie algebra of traceless diagonal matrices. In this case
$\lambda\in \h^*$ can be written as  
$\lambda=(\lambda_1,...,\lambda_n)$, where $\lambda_i\in \mathbb C$. 
Let $v_a$, $a=1,...,n$ be the standard basis of $V$. 
Let $E_{ab}$ be the matrix units 
given by $E_{ab}v_c=\delta_{bc}v_a$. 

We will now give a few examples of quantum dynamical R-matrices.
The general form of the R-matrices will be 
\begin{equation}\label{genform}
R=\sum_a E_{aa}\otimes E_{aa}+\sum_{a\ne b}\alpha_{ab}E_{aa}\otimes E_{bb}
+\sum_{a\ne b}\beta_{ab}E_{ab}\otimes E_{ba},
\end{equation}
where $\alpha_{ab}$, $\beta_{ab}$ are functions which will be given 
explicitly in each example. 

{\bf Example 1.} The basic rational solution. Let
$\beta_{ab}=\frac{1}{\lambda_b-\lambda_a}$,
$\alpha_{ab}=1+\beta_{ab}$. 
Then $R(\lambda)$ is a dynamical R-matrix (with step $1$). 

{\bf Example 2.} The basic trigonometric solution. Let
$\beta_{ab}=\frac{q-1}{q^{\lambda_b-\lambda_a}-1}$, 
$\alpha_{ab}=q+\beta_{ab}$.  
Then $R(\lambda)$ is a dynamical R-matrix (also with step $1$).

{\bf Example 3.} The basic elliptic solution with spectral parameter
(Felder's solution).
Let $\theta(u,\tau)=-\sum_{j\in \mathbb Z+1/2}e^{\pi i(j^2\tau+2j(u+1/2))}$
be the standard theta-function (for brevity will write it as $\theta(u)$). Let 
$\beta_{ab}=\frac{\theta(u-\lambda_b+\lambda_a)\theta(\gamma)}
{\theta(u-\gamma)\theta(\lambda_b-\lambda_a)}$,
$\alpha_{ab}=\frac{\theta(\lambda_a-\lambda_b+\gamma)\theta(u)}{\theta(\lambda_a-
\lambda_b)\theta(u-\gamma)}$. These functions can be viewed as 
functions of $z=e^{2\pi iu}$. 
Then $R(u,\lambda)$ is a quantum dynamical R-matrix 
with spectral parameter and step $\gamma$. 
One may also define the basic trigonometric and basic rational solution 
with spectral parameter, which differ from the elliptic solution 
by replacement of $\theta(x)$ by ${\rm sin}(x)$ and 
$x$, respectively.

{\bf Remark 1.} In examples 1 and 2, the dynamical R-matrix
satisfies the {\it Hecke condition} 
$(PR-1)(PR+q)=0$, with $q=1$ in example 1 (where $P$ 
is the permutation on $V\otimes V$),
and in example 3 the unitarity condition 
$R(u,\lambda)R^{21}(-u,\lambda)=1$. 

{\bf Remark 2.} The basic trigonometric solution 
degenerates into the basic rational solution as $q\to 1$. 
Also, the basic elliptic solution with spectral parameter 
can be degenerated into the basic trigonometric and rational 
solutions with spectral parameter by renormalizing variables and
sending one, respectively both periods 
of theta functions to infinity (see \cite{EV2}).
Furthermore, the basic trigonometric 
and rational solutions with spectral parameter can, in turn, 
be degenerated 
into the solutions of Examples 1,2, by sending the spectral 
parameter to infinity.   
Thus, in essense, all the examples we gave are obtained from 
Felder's solution. 

\subsection{The tensor category of representations associated to a quantum 
dynamical R-matrix}

Let $R$ be a quantum dynamical R-matrix with spectral parameter.
According to \cite{Fe}, 
a representation of $R$ is a semisimple $\h$-module $W$ 
and an invertible meromorphic function 
$L=L_W:\mathbb C\times \h^*\to {\rm End}_\h(V\otimes W)$, such that 
\begin{equation}
\label{dqg}
R^{12}(u_{12},\lambda-h^3)L^{13}(u_{13},\lambda)L^{23}(u_{23},\lambda-h^1)=
L^{23}(u_{23},\lambda)L^{13}(u_{13},\lambda-h^2)R^{12}(u_{12},\lambda).
\end{equation}
(In the case of step $\gamma$, $h^i$ should be replaced by $\gamma h^i$). 

For example: 
$(\mathbb C,1)$ (trivial representation) and $(V,R)$ (vector representation). 

A morphism $f: (W,L_W)\to (W',L_{W'})$ is 
a meromorphic function $f: \h^*\to {\rm End}_\h(W)$ 
such that $(1\otimes f(\lambda))L_W(u,\lambda)=
L_{W'}(u,\lambda)(1\otimes f(\lambda-\gamma h^1))$. 
With this definition, representations form an (additive) category
${\rm Rep}(R)$.
Moreover, it is a tensor category \cite{Fe}: 
given $(W,L_W)$ and $(U,L_U)$, one can form 
the tensor product representation $(W\otimes U,L_{W\otimes U})$, 
where $L_{W\otimes U}(u,\lambda)=L_W^{12}(u,\lambda- \gamma h^3)L_U^{13}
(u,\lambda)$; tensor product of morphisms is defined by 
$(f\otimes g)(\lambda)=f(\lambda-\gamma h^2)\otimes g(\lambda)$.  

In absence of spectral parameter, one should use the same definitions 
without $u$. 

\subsection{Gauge transformations and classification}

There exists a group of rather trivial transformations acting 
on quantum dynamical R-matrices with step $\gamma$. 
They are called {\it gauge transformations}.
If $\h$ and $V$ are as in the previous section, then 
gauge transformations (for dynamical R-matrices without spectral parameter) 
are:

1. Twist by a closed multiplicative 2-form $\phi$:
$\alpha_{ab}\to \alpha_{ab}\phi_{ab}$, 
where $\phi_{ab}=\phi_{ba}^{-1}$, and 
$\phi_{ab}(\lambda)\phi_{bc}(\lambda)\phi_{ca}(\lambda)=
\phi_{ab}(\lambda-\gamma\omega_c)\phi_{bc}(\lambda-\gamma\omega_a)
\phi_{ca}(\lambda-
\gamma\omega_b)$ ($\omega_i={\rm weight}(v_i)$).

2. Permutation of indices $a=1,...,n$;
$\lambda\to \lambda-\nu$;
$R\to cR$. 

In presence of the spectral parameter $u$, 
the constant $c$ in 2 
may depend on $u$, and there are the following additional 
gauge transformations:

3. Given meromorphic $\psi: \h^*\to \mathbb C$, 
$\alpha_{ab}\to e^{u(\psi(\lambda)-2\psi(\lambda-\gamma\omega_a)
+\psi(\lambda-\gamma\omega_a-\gamma\omega_b))}\alpha_{ab}$,\newline 
$\beta_{ab}\to e^{u(\psi(\lambda)-\psi(\lambda-\gamma\omega_a)-
\psi(\lambda-\gamma\omega_b)+\psi(\lambda-\gamma\omega_a-\gamma\omega_b))}
\beta_{ab}$.

4. Multiplication of $u$ by a constant. 

\begin{theorem}\cite{EV2} \label{class1}
Any quantum dynamical $R$-matrix for $\h,V$ satisfying the Hecke 
condition with $q=1$ (respectively, $q\ne 1$) 
is a gauge transformation of the basic rational (respectively, trigonometric) 
solution, or a limit of such $R$-matrices. 
\end{theorem}

In the spectral parameter case, a similar result is known only 
under rather serious restrictions (see \cite{EV2}, Theorem 2.5). 
For more complicated configurations $(\h,V)$, 
classification is not available. 

{\bf Remark.} Gauge transformations 2-4 do not affect the representation 
category of the R-matrix. Gauge transformation 1 does not affect the category 
if the closed form $\phi$ is exact: $\phi=d\xi$, i.e. 
$\phi_{ab}(\lambda)=\xi_a(\lambda)\xi_b(\lambda-\gamma\omega_a)
\xi_a(\lambda-\gamma\omega_b)^{-1}\xi_b(\lambda)^{-1}$, where 
$\xi_a(\lambda)$ is a collection of meromorphic functions.  
We note that this is a very mild condition: 
for example, if $\gamma$ is a formal parameter and we work with 
analytic functions of $\lambda$ in a simply connected domain, then 
a multiplicative 2-form is closed iff it is exact (``multiplicative 
Poincare lemma'').

\subsection{Dynamical quantum groups}

Equation \ref{dqg}
may be regarded as a set of defining relations for an associative algebra
$A_R$ (see \cite{EV2} for precise definitions).
This algebra is a dynamical analogue of the quantum group attached to 
an R-matrix defined in \cite{FRT}, and representations of $R$ are 
an appropriate class of representations of this algebra.
The algebra $A_R$ is called the dynamical quantum group 
attached to $R$. If $R$ is the basic elliptic solution then $A_R$ is 
the elliptic quantum group of \cite{Fe}. 
The structure and representations of $A_R$ 
are studied in many papers (e.g. \cite{Fe,EV2,FV1,TV1}).
  
To keep this paper within bounds, we will not discuss $A_R$ in detail. 
However, let us mention (\cite{EV2})
that $A_R$ is a bialgebroid with base $\h^*$. 
This corresponds to the fact that the category ${\rm Rep}(R)$ is a 
tensor category. Moreover, if $R$ satisfies an additional rigidity assumption 
(valid for example for the basic rational and trigonometric solutions)
then the category ${\rm Rep}(R)$ has duality, and $A_R$ is a Hopf algebroid, 
or a quantum groupoid (i.e. it has an antipode). 

{\bf Remark.} For a general theory of bialgebroids and Hopf algebroids
the reader is referred to \cite{Lu}. However, let us mention that 
bialgebroids with base $X$ correspond to pairs (tensor category,  
tensor functor to $O(X)$-bimodules), similarly to how 
bialgebras correspond to pairs (tensor category, tensor functor to 
vector spaces) (i.e. via Tannakian formalism). 

\subsection{The classical dynamical Yang-Baxter equation}
Recall that if \linebreak $R=1-\hbar r+O(\hbar^2)$ is a solution 
QYBE, then $r$ satisfies the classical Yang-Baxter equation 
(CYBE), 
$[r^{12},r^{13}]+[r^{12},r^{23}]+[r^{13},r^{23}]=0$,
and that $r$ is called the classical limit of $R$,
while $R$ is called a quantization of $r$. 
Similarly, let $R(\lambda,\hbar)$ be a 
family of solutions of QDYBE with step $\hbar$
given by a series $1-\hbar r(\lambda)+O(\hbar^2)$. 
Then it is easy to show that $r(\lambda)$ satisfies 
the following differential equation, 
called the {\it classical dynamical Yang-Baxter equation} (CDYBE)
\cite{BDF,Fe}:
\begin{equation}
\sum_i\biggl(x_i^{(1)}\frac{\partial r^{23}}{\partial x_i}-
x_i^{(2)}\frac{\partial r^{13}}{\partial x_i}+
x_i^{(3)}\frac{\partial r^{12}}{\partial x_i}\biggr)+
[r^{12},r^{13}]+
[r^{12},r^{23}]+
[r^{13},r^{23}]=0,
\end{equation}
(where $x_i$ is a basis of $\h$). 
The function $r(\lambda)$ is called the classical limit 
of $R(\lambda,\hbar)$, and $R(\lambda,\hbar)$ is called a quantization 
of $r(\lambda)$.  

Define a classical dynamical r-matrix to be a meromorphic function 
$r:\h^*\to \text{End}_\h(V\otimes V)$ satisfying CDYBE. 

\begin{conjecture}
Any classical dynamical r-matrix can be quantized. 
\end{conjecture}

This conjecture is proved in \cite{EK} in the non-dynamical case,
and in \cite{Xu} 
in the dynamical case for skew-symmetric solutions
($r^{21}=-r$) satisfying additional technical assumptions. 
However, the most interesting 
non-skew-symmetric case is still open.

Similarly, the classical limit of QDYBE with spectral parameter
is CDYBE with spectral parameter. It is an equation 
with respect to $r(u,\lambda)$ and differs from the usual CDYBE
by insertion of $u_{ij}$ as an additional argument of $r^{ij}$. 

{\bf Remark 1.} Similarly to CYBE, CDYBE makes sense for functions with 
values in $\g\otimes \g$, where $\g$ is a Lie algebra containing $\h$. 

{\bf Remark 2.} The classical limit of the notion of a quantum groupoid 
is the notion of a Poisson groupoid, due to Weinstein. By definition, 
a Poisson groupoid is a groupoid $G$ which is also a Poisson 
manifold, such that the graph of the multiplication is 
coisotropic in $G\times G\times \bar G$ (where $\bar G$ is $G$ with 
reversed sign of Poisson bracket). Such a groupoid 
can be attached (\cite{EV1})
to a classical dynamical r-matrix $r:\h^*\to \g\otimes \g$,
such that $r^{21}+r$ is constant and invariant (i.e. $r$ is a 
``dynamical quasitriangular structure'' on $\g$). This is the 
classical limit of the assignment of a quantum groupoid to 
a quantum dynamical R-matrix (\cite{EV2}). 

\subsection{Examples of solutions of CDYBE}
We will give examples of solutions of CDYBE in the case when $\g$ is a finite 
dimensional simple Lie algebra, and $\h$ is its Cartan subalgebra. 
We fix an invariant inner product on $\g$. 
It restricts to a nondegenerate inner product on $\h$. 
Using this inner product, we identify 
$\h^*$ with $\h$ ($\lambda\in \h^*\mapsto \bar\lambda\in \h$), which yields 
an inner product on $\h^*$. 
The normalization of the inner product 
is chosen so that short roots have squared length $2$.  
Let $x_i$ be an orthonormal basis of $\h$,and 
let $e_\alpha$, $e_{-\alpha}$ denote positive (respectively, negative) root 
elements of $\g$, such that $(e_\alpha,e_{-\alpha})=1$. 

{\bf Example 1.} The basic rational solution
$r(\lambda)=\sum_{\alpha>0}\frac{e_\alpha\wedge e_{-\alpha}}{(\lambda,\alpha)}$.

{\bf Example 2.} The basic trigonometric solution 
$r(\lambda)={\Omega\over 2}+\sum_{\alpha>0}{1\over 2} 
{\rm cotanh}({(\lambda,\alpha)\over 2})$,  
where $\Omega\in S^2\g$ is the inverse element to the inner product on $\g$. 

{\bf Example 3.} The basic elliptic solution with spectral parameter
(Felder's solution) 
$r(u,\lambda)=\frac{\theta'(u)}{\theta(u)}\sum_i x_i\otimes x_i+\sum_{\alpha}
\frac{\theta(u+(\lambda,\alpha))\theta'(0)}{\theta((\lambda,\alpha))
\theta(u)}e_\alpha\otimes e_{-\alpha}$.  

{\bf Remark 1.} One says that a classical 
dynamical r-matrix $r$ has coupling constant $\varepsilon$ 
if $r+r^{21}=\varepsilon \Omega$. If $r$ is with spectral parameter, 
one says that it has 
coupling constant $\varepsilon$ if $r(u,\lambda)+r^{21}(-u,\lambda)=0$, 
and $r(u,\lambda)$ has a simple pole at $u=0$ with residue $\varepsilon
\Omega$ (these are analogs of the 
Hecke and unitarity conditions in the quantum case). 
With these definitions, the basic rational solution 
has coupling constant $0$, while the trigonometric 
and elliptic solutions have coupling constant $1$. 

{\bf Remark 2.} As in the quantum case, Example 3 can be degenerated 
into its trigonometric and rational versions 
where $\theta(x)$ is replaced by $\sin(x)$ and $x$, respectively, 
and Examples 1 and 2 can be 
obtained from Example 3 by a limit. 

{\bf Remark 3.} The classical limit of 
the basic rational and trigonometric solutions of QDYBE
(modified by $\lambda\to \lambda/\hbar$)
is the basic rational, respectively trigoinometric, solutions of CDYBE
for $\g={\mathfrak {sl}}_n$ (in the trigonometric case 
we should set $q=e^{-\hbar/2}$). 
The same is true for the basic elliptic solution 
(with $\gamma=\hbar$). 

{\bf Remark 4.} These examples make sense for any reductive Lie algebra 
$\g$. 

\subsection{Gauge transformations and classification of solutions 
for CDYBE}
It is clear from the above that it is interesting to classify
solutions of CDYBE. As in the quantum case, it should be done up to 
gauge transformations. These transformations are classical analogs 
of the gauge transformations in the quantum case, and look as follows.

1. $r\to r+\omega$, where 
$\omega=\sum_{i,j}C_{ij}(\lambda)x_i\wedge x_j$ is 
a meromorphic closed differential 2-form on $\h^*$. 

2. $r(u,\lambda)\to ar(a\lambda-\nu)$; Weyl group action. 

In the case of spectral parameter, there are additional  
transformations:

3. $u\to bu$. 

4. Let $r=\sum_{i,j} S_{ij}x_i\otimes x_j+\sum_{\alpha}\phi_\alpha e_\alpha 
\otimes e_{-\alpha}$. The transformation is  
$S_{ij}\to S_{ij}+u\frac{\partial^2\psi}{\partial x_i\partial x_j}$, 
$\phi_\alpha\to \phi_\alpha e^{u\partial_\alpha \psi}$, 
where $\psi$ is a function on $\h^*$ with meromorphic $d\psi$. 

\begin{theorem} \cite{EV1}
(i) Any classical dynamical r-matrix with zero coupling constant 
is a gauge transformation of the basic rational solution for 
a reductive subalgebra of $\g$ containing $\h$, or its limiting case.

(ii) Any classical dynamical r-matrix with nonzero coupling constant 
is a gauge transformation of the basic trigonometric solution for 
$\g$, or its limiting case.
 
(iii) Any classical dynamical r-matrix with spectral parameter
and nonzero coupling constant 
is a gauge transformation of the basic elliptic solution for 
$\g$, or its limiting case.
\end{theorem}

{\bf Remark.} One may also classify dynamical r-matrices 
with nonzero coupling constant (without 
spectral parameter) defined on 
${\mathfrak l}^*$ for a Lie subalgebra ${\mathfrak l}
\subset \h$, on which the inner product is nondegenerate (\cite{Sch}). 
Up to gauge transformations 
they are classified by generalized Belavin-Drinfeld triples, 
i.e. triples $(\Gamma_1,\Gamma_2,T)$, where $\Gamma_i$ are 
subdiagrams of the Dynkin diagram $\Gamma$ of $\g$, and 
$T: \Gamma_1\to \Gamma_2$ is a bijection
preserving the inner product of simple roots
(so this classification is a dynamical analog of the 
Belavin-Drinfeld classification of r-matrices on simple Lie algebras, and
the classification 
of \cite{EV1} is the special case $\Gamma_1=\Gamma_2=\Gamma$, $T=1$).
The formula for a classical dynamical r-matrix corresponding to 
such a triple given in \cite{Sch} works for any Kac-Moody algebra, and
in the case of affine Lie algebras yields classical dynamical r-matrices 
with spectral parameter (however, the classification is this case 
has not been worked out). Explicit quantization of the dynamical r-matrices 
from \cite{Sch} (for any Kac-Moody algebra) is given in \cite{ESS}.  

\section{The fusion and exchange construction}
Unlike QYBE, interesting solutions of QDYBE may be obtained already from 
classical representation theory 
of Lie algebras. This can be done through the fusion 
and exchange construction, see \cite{Fa,EV3}.  

\subsection{Intertwining operators} 
Let $\g$ be a simple finite dimensional Lie algebra over $\mathbb C$,
with polar decomposition $\g=\n_+\oplus\h\oplus \n_-$.
For any $\g$-module $V$, we write $V[\nu]$ 
for the weight subspace of $V$ of weight $\nu\in \h^*$. 
Let $M_\lambda$ denote the Verma module over $\g$ with highest weight 
$\lambda\in \h^*$, $x_\lambda$ its highest weight vector, 
and $x_\lambda^*$ the lowest weight vector of the dual module.   
Let $V$ be a finite dimensional representation of $\g$. 
Consider an intertwining operator $\Phi: M_\lambda\to M_\mu\otimes V$. 
The vector $x_\mu^*(\Phi x_\lambda)\in V[\lambda-\mu]$ 
is called the expectation value
of $\Phi$, 
and denoted $\langle \Phi\rangle$. 

\begin{lemma}\label{isom} 
If $M_\mu$ is irreducible (i.e. for generic $\mu$),
the map $\Phi\to \langle \Phi\rangle$ is an isomorphism 
${\rm Hom}_\g(M_{\mu+\nu},M_\mu\otimes V)\to V[\nu]$.  
\end{lemma}

Lemma \ref{isom} allows one to define for any 
$v\in V[\nu]$ (and generic $\lambda$)
the intertwining operator $\Phi_\lambda^v: M_\lambda\to 
M_{\lambda-\nu}\otimes V$, such that $\langle \Phi_\lambda^v\rangle=
v$.

\subsection{The fusion and exchange operators}
Now let $V,W$ be finite dimensional $\g$-modules, and 
$v\in V,w\in W$ homogeneous vectors, of weights ${\rm wt}(v),{\rm wt}(w)$. 
Consider the composition of two intertwining operators 
$$\Phi^{w,v}_\lambda:=
(\Phi^w_{\lambda-{\rm wt}(v)}\otimes 1)\Phi^v_{\lambda}:
M_\lambda\to M_{\lambda-{\rm wt}(v)-{\rm wt}(w)}\otimes W\otimes V.$$ 
The expectation value of this composition, 
$\langle \Phi^{w,v}_\lambda\rangle$, is a bilinear function of $w$ and $v$. 
Therefore, there exists a linear operator 
$J_{WV}(\lambda)\in {\rm End}(W\otimes V)$ 
(of weight zero, i.e. commuting with $\h$), such that 
$\langle \Phi^{w,v}_\lambda\rangle=J_{WV}(\lambda)(w\otimes v)$. 
In other words, we have 
$(\Phi^w_{\lambda-{\rm wt}(v)}\otimes 1)\Phi^v_{\lambda}=
\Phi_\lambda^{J_{WV}(\lambda)(w\otimes v)}$. 
The operator $J_{WV}(\lambda)$ is called the {\it fusion operator}
(because it tells us how to ``fuse'' two intertwining operators).

The fusion operator has a number of interesting properties, which we 
discuss below. In particular, it is lower triangular, i.e. has 
the form $J=1+N$, where $N$ is a sum of terms  
which have strictly positive weights in the second component. Consequently, 
$N$ is nilpotent, and $J$ is invertible. 

Define also the {\it exchange operator} 
$R_{VW}(\lambda):=J_{VW}^{-1}(\lambda)J_{WV}^{21}(\lambda):
V\otimes W\to V\otimes W$. This operator tells us how to 
exchange the order of two intertwining operators, in the sense that 
if $R_{WV}(\lambda)(w\otimes v)=\sum_i w_i\otimes v_i$ (where 
$w_i,v_i$ are homogeneous),
then $\Phi_\lambda^{w,v}=P\sum_i \Phi_\lambda^{v_i,w_i}$
(where $P$ permutes $V$ and $W$).  

\begin{theorem}\cite{EV3}\label{qdybe}
$R_{VV}(\lambda)$ is a solution of QDYBE. 
\end{theorem}

\subsection{Fusion and exchange operators for quantum groups}
The fusion and exchange constructions generalize without significant changes 
to the case when the Lie algebra $\g$ is replaced by the quantum group 
$U_q(\g)$, where $q$ is not zero or a root of unity. The only change 
that needs to be made is in the definition of the exchange operator: 
namely, one sets $R(\lambda)=J_{VW}^{-1}(\lambda)
\mathcal R^{21}J_{WV}^{21}(\lambda)$, where $\mathcal R$ 
is the universal $R$-matrix of $U_q(\g)$. This is because 
when changing the order of intertwining operators, we must 
change the order of tensor product of representations $V\otimes W$,
which in the quantum case is done by means of the R-matrix.  

\begin{example}\label{exch}
Let $\g={\mathfrak{sl}}_n$, and $V$ be the vector representation 
of $U_q(\g)$.
Then the exchange operator has the form $R=q^{1-1/n}\tilde R$, where
$\tilde R$ is given by (\ref{genform}), with 
$\beta_{ab}=\frac{q^{-2}-1}{q^{2(\lambda_a-\lambda_b-a+b)}-1}$, 
$\alpha_{ab}=q^{-1}$ for $a<b$, and 
$\alpha_{ab}={(q^{2(\lambda_b-\lambda_a +a-b)}-q^{-2})
(q^{2(\lambda_b-\lambda_a +a-b)}-q^{2})
\over
q(q^{2(\lambda_b-\lambda_a +a-b)}-1)^2}$ if $a>b$. 
The exchange operator for the vector representation of $\g$
is obtained by passing to the limit $q\to 1$; i.e. it is 
 given by (\ref{genform}), with 
$\beta_{ab}=\frac{1}{\lambda_b-\lambda_a-b+a}$, 
$\alpha_{ab}=1$ for $a<b$, and 
$\alpha_{ab}={(\lambda_b-\lambda_a +a-b-1)
(\lambda_b-\lambda_a +a-b+1)
\over
(\lambda_b-\lambda_a +a-b)^2}$ if $a>b$. 
It is easy to see 
that these exchange operators are gauge equivalent to the basic 
rational and trigonometric solutions of QDYBE, respectively.
\end{example}

\subsection{The ABRR equation}
The fusion operator is not only a tool to define the exchange operator
satisfying QDYBE, but is an interesting object by itself, which deserves
a separate study; so we will briefly discuss its properties. 

Let $\rho$ be the half-sum of 
positive roots of $\g$. Let $\Theta(\lambda)\in U(\h)$ be given by 
$\Theta(\lambda)=\bar\lambda+\bar\rho-\frac{1}{2}\sum x_i^2$.
Then $\Theta(\lambda)$ defines an operator in any 
$U_q(\g)$-module with weight decomposition. Let $\mathcal R_0=
\mathcal R q^{-\sum x_i\otimes x_i}$ be the unipotent part
of the universal $R$-matrix.  

\begin{theorem}\label{ABRR}
(ABRR equation, \cite{ABRR,JKOS}).
For $q\ne 1$, the fusion operator is a unique lower triangular 
zero weight operator, which satisfies the equation:
\begin{equation}
J(\lambda)(1\otimes q^{2\Theta(\lambda)})=
\mathcal R_0^{21}(1\otimes q^{2\Theta(\lambda)})
J(\lambda).  
\end{equation}

For $q=1$, the fusion operator satisfies the classical 
limit of this equation: 
\begin{equation}
[J(\lambda),1\otimes \Theta(\lambda)]=(\sum_{\alpha>0}e_{-\alpha}\otimes 
e_\alpha)J(\lambda),
\end{equation}
\end{theorem}

(Here for brevity we have dropped the subscripts $W$ and $V$,
with the understanding that both sides are operators 
on $W\otimes V$).

\subsection{The universal fusion operator}
Using the ABRR equation, we can define the 
universal fusion operator, living in a completion of $U_q(\g)^{\otimes 2}$, 
which becomes $J_{WV}(\lambda)$ 
after evaluating in $W\otimes V$. Namely, the universal fusion operator 
$J(\lambda)$ is the unique lower triangular solution 
of the ABRR equation in a completion of $U_q(\g)^{\otimes 2}$.
This solution can be found in the form of a series
$J=\sum_{n\ge 0} J_n$, $J_0=1$, where $J_n\in U_q(\g)\otimes U_q(\g)$ 
has zero weight and its second component has degree $n$ in principal 
gradation; so $J_n$ are computed recursively. 

This allows one to compute the universal fusion operator quite explicitly. 
For example, if $q=1$ and $\g={\mathfrak sl}_2$, then 
the universal fusion operator is given by the formula
$J(\lambda)=\sum_{n\ge 0}
\frac{(-1)^n}{n!}f^n \otimes (\lambda-h+n+1)^{-1} \ldots 
(\lambda-h +2n)^{-1}e^n$.

\subsection{The dynamical twist equation}
Another important property of the fusion operator 
is the dynamical twist equation (which 
is a dynamical analog of the equation for a Drinfeld 
twist in a Hopf algebra). 

\begin{theorem}\label{dytwist}
The universal fusion operator $J(\lambda)$ satisfies 
the dynamical twist equation 
$J^{12,3}(\lambda)J^{1,2}(\lambda-h^{(3)})=J^{1,23}(\lambda)J^{2,3}(\lambda)$.
\end{theorem}

Here the superscripts of $J$ stand for components where the first 
and second component of $J$ acts; for example $J^{1,23}$ means 
$(1\otimes \Delta)(J)$, and $J^{1,2}$ means $J\otimes 1$. 

\subsection{The fusion operator for affine algebras}
The fusion and exchange construction can be generalized to the case 
when the simple Lie algebra $\g$ is replaced by any Kac-Moody Lie algebra. 
This generalization is especially interesting if $\g$ is replaced 
with the affine Lie algebra $\hat\g$, and $V,W$ are finite dimensional 
representations of $U_q(\hat\g)$ (where $q$ is allowed to be $1$). In this case, 
for each $z\in \mathbb C^*$ we have an outer automorphism $D_z: U_q(\hat\g)\to 
U_q(\hat\g)$, which preserves 
the Chevalley 
generators $q^h$ and $e_i,f_i$, $i>0$, of $U_q(\hat\g)$,
while multiplying $e_0$ by $z$ 
and $f_0$ by $z^{-1}$. Define the shifted representation $V(z)$ by 
$\pi_{V(z)}(a)=\pi_V(D_z(a))$. Let $\hat\lambda=(\lambda,k)$ be a weight for 
$\hat\g$ ($k$ is the level). Then similarly to the finite dimensional
case one can define the intertwining operator 
$\Phi_{\lambda,k}^v(z): M_{\hat\lambda}\to M_{\hat\lambda-{\rm wt}(v)}
\hat\otimes V(z)$ (to a completed tensor product);
it can be written as an infinite series $\sum_{n\in \mathbb Z}
\Phi_{\lambda,k}^v[n]z^{-n}$, 
where $\Phi_{\lambda,k}^v[n]:
M_{\hat\lambda}\to M_{\hat\lambda-{\rm wt}(v)}
\otimes V$ are linear operators. Furthermore, 
the expectation value of the composition 
$\langle (\Phi^w_{\lambda-{\rm wt}(v),k}(z_1)\otimes 1)\Phi^v_{\lambda,k}(z_2)
\rangle$ is defined as a Taylor series in $z=z_2/z_1$. Thus, one can define
the fusion operator $J_{WV}(z)\in {\rm End}_\h(W\otimes V)[[z]]$
such that this expectation value is equal to $J_{WV}(z)(w\otimes v)$. 

\subsection{Fusion operator and correlators in the WZW model}
One may show (see \cite{FR,EFK}) 
that the series $J_{WV}(z)$ is convergent 
in some neighborhood of $0$ to a holomorphic function. 
This function has a physical interpretation. 
Namely, if $q=1$, the operators $\Phi_{\lambda,k}^v(z)$ are, 
essentially, vertex
operators (primary fields) for the Wess-Zumino-Witten conformal field theory, 
and the function $J_{WV}(z)(w\otimes v)$ 
is a 2-point correlation function of vertex operators. 
If $q\ne 1$, this function is a 2-point correlation 
function of q-vertex operators, and has a similar interpretation in terms
of statistical mechanics. 

\subsection{The ABRR equation in the affine case as the KZ (q-KZ) equation}
The ABRR equation has a straightforward generalization to the affine case, 
which also has a physical interpretation. Namely, in the case 
$q=1$ it coincides with the (trigonometric) Knizhnik-Zamolodchikov (KZ) 
equation 
for the 2-point correlation function, while for $q\ne 1$
it coincides with the quantized KZ equation 
for the 2-point function of q-vertex operators, derived in 
\cite{FR} (see also \cite{EFK}). 

One may also define (for any Kac-Moody algebra) the {\it multicomponent 
universal fusion operator}
$J^{1...N}(\lambda)=J^{1,2...N}(\lambda)...J^{2,3...N}(\lambda)
J^{N-1,N}(\lambda)$. 
It satisfies a multicomponent version of the ABRR equation.
In the affine case, $J^{1...N}(\lambda)$ is interpreted as 
the N-point correlation function for vertex (respectively, q-vertex)
operators, while the multicomponent ABRR equation is interpreted as 
the KZ (respectively, qKZ) equation for this function. 
See \cite{EV5} for details. 

\subsection{The exchange operator 
for affine algebras and monodromy of KZ (q-KZ) equations} 
The generalization of the exchange operator
to the affine case is rather tricky, and there is a serious difference 
between the classical ($q=1$), and quantum ($q\ne 1$) case. 
The naive definition would be $R_{VW}(u,\hat\lambda)=
J_{VW}(z,\hat\lambda)^{-1}\mathcal R^{21}|_{V\otimes W(z)}J_{WV}^{21}
(z^{-1},\hat\lambda)$ (where $z=e^{2\pi iu}$).
However, this definition does not immediately make sense, since 
we are multiplying infinite Taylor series in $z$ by infinite 
Taylor series in $z^{-1}$. To make sense of such product, 
consider the cases $q=1$ and $q\ne 1$ separately. 

If $q=1$, the functions $J_{VW}(z,\hat\lambda)$ and 
$J_{WV}^{21}(z^{-1},\hat\lambda)$ are 
both solutions of the Knizhnik-Zamolodchikov differential 
equation, one regular near $0$ and the other near $\infty$.
The equation is nonsingular outside $0,\infty$, and $1$.  
Thus, the series $J_{VW}(z,\hat\lambda)$ is convergent for $|z|<1$. 
On the other hand, for $0<z<1$ define 
the function $A^{\pm}(J_{WV}^{21}(z^{-1},\hat\lambda))$
to be the analytic continuation of $J_{WV}^{21}(z^{-1},\hat\lambda)$
from the region $z>1$ along a curve passing $1$ from above (for +) and below
(for -). Then the function 
$R_{VW}^\pm(u,\hat\lambda):=
J_{VW}(z,\hat\lambda)^{-1}A^\pm(J_{WV}^{21}
(z^{-1},\hat\lambda))$  
is of zero weight, and satisfies the QDYBE with spectral parameter
(for $V=W$).
The operator $R_{VW}^\pm(u,\hat\lambda)$ 
is the appropriate analog of the exchange operator
(depending on a choice of sign). 

If $q\ne 1$, the functions $J_{VW}(z,\hat\lambda)$ and 
$\mathcal R^{21}|_{V\otimes W(z)}J_{WV}^{21}(z^{-1},\hat\lambda)$ are 
both solutions of the quantized 
Knizhnik-Zamolodchikov 
equation, which is a difference equation with multiplicative step 
$p=q^{-2m(k+g)}$, where $m$ is the ratio of squared norms of long and short 
roots, $k$ the level of $\hat\lambda$, and $g$ the dual Coxeter number 
of $\g$. Therefore, if $|p|\ne 1$, these functions admit a meromorphic 
continuation to the whole $\mathbb C^*$ (unlike the $q=1$ case, they are now 
single-valued), and the naive formula 
for $R_{VW}(u,\hat\lambda)$ makes sense. As in the $q=1$ case, this 
function is of zero weight and satisfies the QDYBE with spectral parameter
(for $V=W$); it is the appropriate generalization of the exchange operator
(see \cite{EFK} for details). 

The operators $R_{VW}^\pm(z,\hat\lambda)$, essentially,
represent the monodromy of the KZ differential equation.
In particular, $R^\pm_{VW}$  is ``almost constant'' in $u$:
its matrix elements in a homogeneous basis under 
$\h$ are powers of $e^{2\pi iu}$ (in fact, 
$R_{VV}^\pm(u,\hat\lambda)$ 
is gauge equivalent, in an appropriate sense, to a 
solution of QDYBE without spectral parameter). 
Similarly,
the operator $R_{VW}(u,\hat\lambda)$ for $q\ne 1$ represents 
the q-monodromy (Birkhoff's connection matrix) 
of the q-KZ equation. In particular, the matrix elements 
of $R_{VW}(u,\hat\lambda)$ 
are quasiperiodic in $u$ with period $\tau=\log p/2\pi i$. 
Since they are also periodic with period $1$ and meromorphic, 
they can be expressed rationally via elliptic theta-functions. 
For $\g={\frak{sl}}_n$, they were computed in \cite{TV2}.

\begin{example}\label{exaff}
 Let $\g={\mathfrak {sl}}_n$, and $V$ the vector representation
of $U_q(\hat\g)$. 
If $q\ne 1$, the exchange operator $R_{VV}(u,\hat\lambda)$ 
is a solution of QDYBE, gauge equivalent to the basic elliptic 
solution with spectral parameter (see \cite{Mo}
and references therein, 
 and also \cite{FR,EFK}).
The gauge transformation involves an exact multiplicative 2-form, 
which expresses via q-Gamma functions with ${\rm q}=p$. 
Similarly, if $q=1$, the exchange operators $R_{VV}^\pm(u,\hat\lambda)$
are gauge equivalent to the basic trigonometric solution without 
spectral parameter, with $q=e^{\pi i/m(k+g)}$; the gauge transformation 
involves an exact multiplicative 2-form expressing via classical 
Gamma-functions. This is obtained by 
sending $q$ to $1$ in the result of \cite{Mo}.
\end{example}

{\bf Remark.} Note that the limit $q\to 1$ in this setting is rather subtle. 
Indeed, for $q\ne 1$ the function $J_{VW}(u,\lambda)$ has an infinite 
sequence of poles in the z-plane, which becomes denser as 
$q$ approaches $1$ and eventually degenerates into a branch 
cut; i.e. this single valued meromorphic function becomes multivalued 
in the limit.

\subsection{The quantum Kazhdan-Lusztig functor}
Let $\g={\mathfrak{sl}}_n$, $q\ne 1$. 
Example \ref{exaff} allows one to construct a tensor functor 
from the category ${\rm Rep}_f(U_q(\hat\g))$ 
of finite dimensional $U_q(\hat\g)$-modules, 
to the category ${\rm Rep}_f(R)$ of finite dimensional 
representations of the basic elliptic solution $R$ of QDYBE 
with spectral parameter (i.e. to the category of finite dimensional 
representations of Felder's elliptic quantum group). 
Namely, let $V$ be the vector representation, 
and for any finite dimensional representation $W$ of $U_q(\hat\g)$, 
set $L_W=R_{VW}$. Then $(W,L_W)$ is a representation of the dynamical 
R-matrix $R_{VV}$. This defines a functor $F: {\rm Rep}_f(U_q(\hat\g))\to   
{\rm Rep}_f(R_{VV})$.  
Moreover, this functor is a tensor functor: the tensor structure 
$F(W)\otimes F(U)\to F(W\otimes U)$ 
is given by the fusion operator $J_{WU}(\hat\lambda)$ (the axiom  
of a tensor structure follows from the dynamical twist equation for $J$).
On the other hand, since $R_{VV}$ and $R$ are gauge equivalent by an exact 
form, their representation categories are equivalent, so 
one may assume that $F$ lands in ${\rm Rep}_f(R)$. 

If the scalars for ${\rm Rep}_f(U_q(\hat\g))$ 
are taken to be the field of periodic functions of $\hat\lambda$
(in particular, $k$ is regarded as a variable), then
$F$ is fully faithful; it can be regarded as a q-analogue of the 
Kazhdan-Lusztig functor  
(see \cite{EM} and references therein).
It generalizes to infinite dimensional representations, and 
allows one to construct elliptic deformations of all evaluation 
representations of $U_q(\hat\g))$ (which was done for finite dimensional 
representations in \cite{TV1}). 
The versions of this functor without spectral parameter,
 from representations of $\g$ ($U_q(\g)$) 
to representations of the basic rational (trigonometric) solution of QDYBE
can be found in \cite{EV3}. 

\section{Traces of intertwining operators and Macdonald functions}
In this section we discuss a connection between dynamical R-matrices and 
certain integrable systems and special functions (in particular, 
Macdonald functions).

\subsection{Trace functions}
Let $V$ be a finite dimensional representation of $U_q(\g)$ ($q\ne 1$), 
such that $V[0]\ne 0$. Recall that for 
any $v\in V[0]$ and generic $\mu$, 
one can define an intertwining operator $\Phi_\mu^v$
such that $\langle \Phi_\mu^v\rangle=v$.
Following \cite{EV4} set
$\Psi^v(\lambda,\mu)={\rm Tr}|_{M_\mu}(\Phi_\mu^v q^{2\bar\lambda})$.
This is an infinite series in the variables $q^{-(\lambda,\alpha_i)}$
(where $\alpha_i$ are the simple roots)
whose coefficients are rational functions of $q^{(\mu,\alpha_i)}$
(times a common factor $q^{2(\lambda,\mu)}$). 
For generic $\mu$ 
this series converges near $0$, and its matrix elements belong to
 $q^{2(\lambda,\mu)}(\mathbb C(q^{(\lambda,\alpha_i)})
\otimes \mathbb C(q^{(\mu,\alpha_i)}))$. 

Let $\Psi_V(\lambda,\mu)$ be the ${\rm End}(V[0])$-valued
function, such that $\Psi_V(\lambda,\mu)v=\Psi^v(\lambda,\mu)$. 
The function $\Psi_V$ has remarkable properties and in a special case is 
closely related to Macdonald functions. To formulate the properties of 
$\Psi_V$, we will consider a renormalized version of this function.
Namely, let $\delta_q(\lambda)$ be the Weyl denominator 
$\prod_{\alpha>0}(q^{(\lambda,\alpha)}-q^{-(\lambda,\alpha)})$. 
Let also $Q(\mu)=\sum S^{-1}(b_i)a_i$, where 
$\sum a_i\otimes b_i=J(\mu)$ is the universal fusion operator
(this is an infinite expression, but it makes sense  
as a linear operator on finite dimensional 
representations; moreover it is of zero weight and invertible). 
Define the {\it trace function} $F_V(\lambda,\mu)=
\delta_q(\lambda)\Psi_V(\lambda,-\mu-\rho)Q(-\mu-\rho)^{-1}$. 

\subsection{Commuting difference operators} 
For any finite dimensional $U_q(\g)$-module $W$,
we define a difference operator $\mathcal D_W$ acting on functions 
on $\h^*$ with values in $V[0]$. Namely, we set 
$(\mathcal D_W f)(\lambda)=\sum_{\nu\in \h^*}{\rm 
Tr}|_W(R_{WV}(-\lambda-\rho))f(\lambda+\nu)$. 
These operators are dynamical analogs of transfer matrices, and were 
introduced in \cite{FV2}. It can be shown that 
${\mathcal D}_{W_1\otimes W_2}={\mathcal D}_{W_1}{\mathcal D}_{W_2}$;
in particular, ${\mathcal D}_W$ commute with each other, and the algebra 
generated by them is the polynomial algebra in 
${\mathcal D}_{W_i}$, where $W_i$ are the fundamental representations
of $U_q(\g)$. 

\subsection{Difference equations for the trace functions}
It turns out that trace functions $F_V(\lambda,\mu)$, regarded as 
functions of $\lambda$,
are common eigenfunctions of 
${\mathcal D}_W$.
\begin{theorem}\cite{EV4}\label{macd} 
One has 
\begin{equation}\label{macd1}
{\mathcal D}_W^{(\lambda)}F_V(\lambda,\mu)=
\chi_W(q^{-2\bar\mu})F_V(\lambda,\mu),
\end{equation}
where $\chi_W(x)={\rm Tr}|_W(x)$ is the character of $W$.
\end{theorem}
In fact, it is easy to deduce from this theorem that 
if $v_i$ is a basis of $V[0]$ then $F_V(\lambda,\mu)v_i$ is a basis 
of solutions of (\ref{macd1}) in the power series space. 
Thus, trace functions allow us to integrate the quantum integrable system 
defined by the commuting operators ${\mathcal D}_{W_i}$. 
\begin{theorem}\label{macd2}\cite{EV4} 
The function $F_V$ is symmetric in $\lambda$ and $\mu$ 
in the following sense: $F_{V^*}(\mu,\lambda)=F_V(\lambda,\mu)^*$. 
\end{theorem}
This symmetry property implies that 
$F_V$ also satisfies ``dual'' difference equations 
with respect to $\mu$:
${\mathcal D}_W^{(\mu)}F_V(\lambda,\mu)^*=
\chi_W(q^{-2\bar\lambda})F_V(\lambda,\mu)^*$.

\subsection{Macdonald functions}
An important special case, worked out in \cite{EKi}, 
is $\g={\mathfrak sl}_n$, and 
$V=L_{kn\omega_1}$, where $\omega_1$ is the first fundamental weight, 
and $k$ a nonnegative integer. In this case, ${\rm dim}V[0]=1$, and 
thus trace functions can be regarded as scalar functions. 
Furthermore, it turns out (\cite{FV2}) that the operators 
${\mathcal D}_W$ can be conjugated (by a certain explicit product) 
to Macdonald's difference operators 
of type $A$, and thus the functions $F_V(\lambda,\mu)$ 
are Macdonald functions (up to multiplication by this product).
One can also obtain Macdonald's polynomials by replacing Verma modules 
$M_\lambda$ with irreducible finite dimensional modules $L_\lambda$;
see \cite{EV4} for details. In this case, Theorem \ref{macd2} 
is the well known 
Macdonald's symmetry identity, and the ``dual'' difference equations are
the recurrence relations for Macdonald's functions. 

{\bf Remark 1.} The dynamical transfer matrices ${\mathcal D}_W$ can be 
constructed not only for the trigonometric 
but also for the elliptc dynamical R-matrix; in the case 
$\g={\mathfrak {sl}}_n$, $V=L_{kn\omega_1}$ this yields the Ruijsenaars 
system, which is an elliptic deformation of the Macdonald system. 

{\bf Remark 2.} If $q=1$, the difference equations of Theorem \ref{macd}
become differential equations, which in the case 
$\g={\mathfrak {sl}}_n$, $V=L_{kn\omega_1}$
reduce to the trigonometric Calogero-Moser system. 
In this limit, the symmetry property 
is destroyed, but the ``dual'' difference equations remain valid, 
now with the exchange operator for $\g$ rather than $U_q(\g)$. 
Thus, both for $q=1$ and $q\ne 1$, common eigenfunctions satisfy 
additional difference equations with respect to eigenvalues -- the so called 
bispectrality property. 

{\bf Remark 3.} Apart from trace $\Psi^v$ of a single intertwining operator 
multiplied by $q^{2\bar\lambda}$, it is useful to consider 
the trace of a product of several such operators. After an appropriate 
renormalization, such multicomponent trace function 
(taking values in ${\rm End}((V_1\otimes...\otimes V_N)[0]$) satisfies 
multicomponent analogs of (\ref{macd1}) and its dual version, 
as well as the symmetry. Furthermore, it satisfies an additional 
quantum Knizhnik-Zamolodchikov-Bernard equation, and its dual version
(see \cite{EV4}). 

{\bf Remark 4.} The theory of trace functions can be generalized 
to the case of affine Lie algebras, with $V$ being 
a finite dimensional representation of $U_q(\hat\g)$. 
In this case, trace functions will be interesting transcendental functions.
In the case $\g={\mathfrak {sl}}_n$, $V=L_{kn\omega_1}$, they are 
analogs of Macdonald functions for affine root systems. 
It is expected that for $\g={\mathfrak {sl}_2}$ they are the 
elliptic hypergeometric functions studied in \cite{FV3}. 
This is known in the trigonometric limit (\cite{EV4})
and for $q=1$. 

{\bf Remark 5.} The theory of trace functions, both finite dimensional 
and affine, can be generalized to the case of any 
generalized Belavin-Drinfeld triple; see \cite{ES1}.

{\bf Remark 6.} Trace functions $F_V(\lambda,\mu)$ are not Weyl group 
invariant. Rather, the diagonal action of the Weyl group multiplies them 
by certain operators, called the {\it dynamical Weyl group} operators
(see \cite{TV3,EV5}). These operators play an important role 
in the theory of dynamical R-matrices and trace functions, 
but are beyond the scope of this paper.

\bibliographystyle{ams-alpha}

\begin{thebibliography}{A}  

\bibitem[ABRR]{ABRR}
 D.Arnaudon, E.Buffenoir, E.Ragoucy, and Ph.Roche,  
\emph{Universal Solutions of quantum dynamical Yang-Baxter equations},
Lett. Math. Phys. \textbf{44} (1998), no. 3, 201-214.

\bibitem [BDF]{BDF} J. Balog, L. Dabrowski, and L. Feh\'er,
Classical $r$-matrix and exchange algebra in WZNW and 
Toda field theories,
Phys. Lett. B, v. 244, issue 2, p.227--234, 1990.

\bibitem[Dr]{Dr}
Drinfeld, V. G. Quantum groups. Proceedings of the 
International Congress of Mathematicians, Vol. 1, 2 (Berkeley, Calif., 1986),
798--820, Amer. Math. Soc., Providence, RI, 1987.

\bibitem[EFK]{EFK} 
Etingof P.,  Frenkel I., Kirillov Jr. A.,
\emph{Lectures on representation theory 
and Knizhnik-Zamolodchikov equations} AMS, (1998).

\bibitem[EK]{EK}
Etingof P., Kazhdan D., \emph{Quantization of Lie bialgebras I},
Selecta Math. \textbf{2} (1996), 1-41.

\bibitem[EKi]{EKi}
 Etingof, P.I., Kirillov, A.A., Jr,  \emph{Macdonald's
polynomials and representations of quantum groups}, 
Math. Res. Let., \textbf{1} (3), (1994) 279-296.

\bibitem[EM]{EM}
Etingof, P., Moura, A.,
On the quantum Kazhdan-Lusztig functor,
math.QA/0203003, 2002.

\bibitem[ES1]{ES1}
Etingof, P.; Schiffmann, O. Twisted traces of quantum intertwiners and quantum dynamical $R$-matrices corresponding to generalized
Belavin-Drinfeld triples. CMP 218 (2001), no. 3, 633-663. 

\bibitem[ES2]{ES2}
Etingof, P.; Schiffmann, O. 
Lectures on the dynamical Yang-baxter equations, math.QA/9908064.

\bibitem[ESS]{ESS}
Etingof P., Schedler T., Schiffmann O., \emph{Explicit quantization of dynamical r-matrices}, J. Amer. Math. Soc., \textbf{13} 595-609, (2000).

\bibitem[EV1]{EV1}
Etingof P., Varchenko A., \emph{Geometry and classification of solutions of the classical dynamical Yang-Baxter equation}, Commun. Math. Phys, \textbf{192} 77-120 (1998).

\bibitem[EV2]{EV2}
Etingof P., Varchenko A., \emph{Solutions of the quantum dynamical Yang-Baxter equation and dynamical quantum groups}, Commun. Math. Phys, \textbf{196} 591-640 (1998).

\bibitem[EV3]{EV3}
Etingof P., Varchenko A., \emph{Exchange dynamical quantum groups},  CMP 205 (1999), no. 1, 19--52. 

\bibitem[EV4]{EV4}
Etingof P., Varchenko A., \emph{Traces of intertwiners for quantum groups
and difference equations, I}, 
 Duke Math. J. 104 (2000), no.
3, 391--432.

\bibitem[EV5]{EV5}
Etingof P., Varchenko A., \emph{Dynamical Weyl groups and applications}, 
math.QA/0011001.

\bibitem[Fa]{Fa}
 Faddeev L.,  \emph{On the exchange matrix 
of the WZNW model}, CMP, \textbf{132} (1990), 131-138.

\bibitem[Fe]{Fe}
 Felder G., \emph{Conformal field theory and integrable
systems associated to elliptic curves},
Proceedings of the International Congress of Mathematicians,
Z\"urich 1994, p.1247--1255.

\bibitem[FR]{FR}
Frenkel I., Reshetikhin N.,
\emph{Quantum affine algebras and holonomic
difference equations},
Commun. Math. Phys. \textbf{146} (1992), 1-60.

\bibitem[FRT]{FRT}
Faddeev, L.D.; 
Reshetikhin, N. Yu.; 
Takhtajan, L. A. Quantization of Lie groups and Lie algebras. Algebraic analysis, Vol. I, 129--139,
Academic Press, Boston, MA, 1988. 


\bibitem[FV1]{FV1} 
Felder G., Varchenko A.,
\emph{On representations of the elliptic quantum group
$E_{\tau,\eta}(sl_2)$}, Commun. Math. Phys. \textbf{181} (1996),
746--762.

\bibitem[FV2]{FV2}
Felder G., Varchenko A., 
\emph{Elliptic quantum groups and Ruijsenaars models},
J. Statist. Phys. \textbf{89} (1997), no. 5-6, 963-980.

\bibitem[FV3]{FV3} 
Felder G., Varchenko A.,
The $q$-deformed Knizhnik-Zamolodchikov-Bernard heat equation. 
CMP 221 (2001),
no. 3, 549--571.

\bibitem[JKOS]{JKOS}
Jimbo, M.; Odake, S.; Konno, H.; Shiraishi, J. Quasi-Hopf twistors for elliptic quantum groups. Transform. Groups 4 (1999), no. 4, 303--327.

\bibitem[Lu]{Lu} 
Lu J.H.,  \emph{Hopf algebroids and quantum groupoids}, 
Inter. J. Math., \textbf{7} (1), (1996), 47-70.

\bibitem[Mo]{Mo} Moura, A.,
Elliptic Dynamical R-Matrices from the Monodromy of the q-Knizhnik-Zamolodchikov Equations for the Standard Representation of Uq(sl(n+1)),
math.RT/0112145.

\bibitem[Sch]{Sch}
Schiffmann O, \emph{On classification of dynamical r-matrices}, MRL, \textbf{5}, 13-30 (1998).

\bibitem[TV1]{TV1}
Tarasov, V.; Varchenko, A. Small elliptic quantum group $e_{\tau,\gamma}
({\mathfrak{sl}}\sb N)$. Mosc. Math. J. 1 (2001), no. 2, 243--286,
303--304.

\bibitem[TV2]{TV2}
Tarasov, V.; Varchenko, A.
Geometry of $q$-hypergeometric 
functions, quantum affine algebras and elliptic quantum groups. Astérisque No.
246 (1997).

\bibitem[TV3]{TV3}
Tarasov, V.; Varchenko, A.
Difference equations compatible with trigonometric KZ differential equations. 
IMRN 2000,
no. 15, 801--829.

\bibitem[Xu]{Xu} Xu, P.,
Triangular dynamical $r$-matrices and quantization. Adv. Math. 166 (2002), no. 1, 1--49.

\end{thebibliography}

\end{document}